\documentclass[a4paper]{article}
\usepackage{amssymb,amsmath,amscd,amsthm}
\usepackage[all,cmtip]{xy}
\usepackage{multirow}

\title{Essential Cohomology of the p-Groups with a Cyclic Subgroup of Index p}
\author{\\Christopher A. Gerig\\
\textit{Cornell University}}
\date{}

\begin{document}
\maketitle

\noindent\textbf{Abstract.}  In this paper we determine explicitly the mod-$p$ essential cohomology ideals of the $p$-groups with a cyclic subgroup of index $p$.

\setcounter{section}{-1}
\section{Introduction}
$\;\;\;\;\;\;$ Let $G$ be a finite group and $k$ the field $\mathbb{Z}_p$.  An element $x\in H^*(G,k)$ is called \textit{essential} if it has trivial restriction to all proper subgroups of $G$; these elements make up the essential ideal $Ess(G)$.  This ideal measures the failure of the set of maximal subgroups to detect $H^*(G,k)$, where a collection $X$ of subgroups of $G$ \textit{detects} $H^*(G,k)$ if the induced map by restrictions $H^*(G,k)\rightarrow \prod_{H\in X}H^*(H;k)$ is injective.  Other important properties and further information may be found in [Ad2, Gr, Ma].\\
\indent In the literature, both theory and computation have been provided on the cohomology rings of groups.  However, only theory has been provided on the topic of Essential Cohomology, with the exceptions of the elementary abelian $p$-groups (see [Ak]) and extraspecial $p$-groups (see [Mi2]).  This paper begins to fill in the gap by providing the essential ideals of the $p$-groups which have a cyclic subgroup of index $p$.\\
\indent The outline of this paper is as follows.  In section 1 we state the classification of $p$-groups with a cyclic subgroup of index $p$, followed by their mod-$p$ cohomology rings.  In section 2 we list a few relevant facts on essential cohomology, and the subsequent sections contain the calculations of the essential ideals.

\section{Preliminaries}
The classification of the $p$-groups with a cyclic subgroup of index $p$ is well-known, and is proved in [Br] using the cohomology theory of extensions.  For convenience we list the theorem here.\\
\\
\textbf{Theorem 1.1.}  \textit{If G is a p-group with a cyclic subgroup of index p, then G is isomorphic to one of the following groups:\\
(A)  Cyclic:\indent $\mathbb{Z}_{p^n}\;,\;\;n\ge 1$\\
\indent\indent $\langle t\;|\;t^{p^n}=1\rangle$\\
(B)  Direct Product:\indent $\mathbb{Z}_{p^n}\times\mathbb{Z}_p\;,\;\;n\ge 1$\\
\indent\indent $\langle t,s\;|\; t^{p^n}=s^p=1,\, st=ts\rangle$\\
(C)  Nonabelian Split Metacyclic:\indent $\mathbb{Z}_{p^n}\rtimes \mathbb{Z}_p\;,\;\;n\ge 2$\\
\indent\indent $\langle t,s\;|\;s^p=t^{p^n}=1,\,sts^{-1}=t^{p^{n-1}+1}\rangle$\\
(D)  Dihedral:\indent $D_{2^n}\;,\;\; n\ge 3$\\
\indent\indent $\langle t,s\;|\; t^{2^{n-1}}=s^2=(ts)^2=1\rangle $\\
(E)  Generalized Quaternion:\indent $Q_{2^n}\;,\;\; n\ge 3$\\
\indent\indent $\langle t,s\;|\; t^{2^{n-2}}=s^2,\, t^{2^{n-1}}=1,\, tst=s\rangle$\\
(F)  Nonabelian Split Metacyclic:\indent $\mathbb{Z}_{2^n}\rtimes \mathbb{Z}_2\;,\;\;n\ge 3$\\
\indent\indent $\langle t,s\;|\;s^2=t^{2^n}=1,\,sts=t^{2^{n-1}-1}\rangle$}\\
\\
\indent The mod-p cohomology rings of these groups are also well known, and will be stated here without proof (classes C+F are given in [Di], classes D+E are given in [Ad1], and class B follows immediately from class A using the K$\ddot{u}$nneth formula).
\\
\\
\indent\textbf{Class A}\indent\indent\indent\indent\indent\indent\indent\indent\indent\textbf{Class B}\\
$H^*(\mathbb{Z}_{p^n},k)\cong \Lambda_k[x]\otimes_k k[y]$\indent\indent\indent$H^*(\mathbb{Z}_{p^n}\times\mathbb{Z}_p,k)\cong\Lambda_k[x_1,x_2]\otimes_k k[y_1,y_2]$\\
\indent where $p^n\ne 2$\indent\indent\indent\indent\indent\indent\indent where $n\ge 1$, $p>2$\\
\indent $|x|=1$, $|y|=2$\indent\indent\indent\indent\indent\indent\indent $|x_1|=|x_2|=1$, $|y_1|=|y_2|=2$\\
$H^*(\mathbb{Z}_{2},k)\cong k[x]$\indent\indent\indent\indent\indent\indent$H^*(\mathbb{Z}_{2^n}\times\mathbb{Z}_2,k)\cong \Lambda_k[x_1]\otimes_k k[x_2,y_1]$\\
\indent $|x|=1$\indent\indent\indent\indent\indent\indent\indent\indent\indent where $n\ge 2$ \\
\indent\indent\indent\indent\indent\indent\indent\indent\indent\indent\indent\indent $|x_1|=|x_2|=1$, $|y_1|=2$\\
\indent\indent\indent\indent\indent\indent\indent\indent\indent\indent\indent$H^*(\mathbb{Z}_2\times\mathbb{Z}_2,k)\cong k[x_1,x_2] $\\
\indent\indent\indent\indent\indent\indent\indent\indent\indent\indent\indent\indent $|x_1|=|x_2|=1$\\
\\
\indent\textbf{Class C}\\
$H^*(\mathbb{Z}_{p^n}\rtimes \mathbb{Z}_p,k)\cong k[a_1,\ldots,a_{p-1},b,y,v,w]/(b^2,v^2,a_ia_j,a_iv,a_iy)$\\
\indent where $n\ge 2$, $p>2$\\
\indent $|a_i|=2i-1$, $|b|=1$, $|y|=2$, $|v|=2p-1$, $|w|=2p$\\
$H^*(\mathbb{Z}_{2^n}\rtimes \mathbb{Z}_2,k)\cong k[a,b,v,w]/(a^2,v^2,av,ab^2)$\\
\indent where $n\ge 2$\\
\indent $|a|=|b|=1$, $|v|=3$, $|w|=4$\\
\\
\indent\textbf{Class D}\\
$H^*(D_{2^n},k)\cong k[x,y,z]/(xy)$\\
\indent where $n\ge 3$\\
\indent $|x|=|y|=1$, $|z|=2$\\
\\
\indent\textbf{Class E}\\
$H^*(Q_8,k)\cong k[x,y,z]/(x^2+xy+y^2,x^2y+xy^2)$\\
\indent $|x|=|y|=1$, $|z|=4$\\
$H^*(Q_{2^n},k)\cong k[x,y,z]/(xy,x^3+y^3)$\\
\indent where $n\ge 4$\\
\indent $|x|=|y|=1$, $|z|=4$\\
\\
\indent\textbf{Class F}\\
$H^*(\mathbb{Z}_{2^n}\rtimes \mathbb{Z}_2,k)\cong k[a,b,y,v,w]/(ay,av,b^2,a^2+ab,v^2+wab+vyb)$\\
\indent where $n\ge 3$\\
\indent $|a|=|b|=1$, $|y|=2$, $|v|=3$, $|w|=4$

\newpage
\section{Facts on Essential Cohomology}
$\;\;\;\;\;\;$ Consider the cohomology ring $H^*(G,k)$ of a finite group $G$ with coefficients $k\equiv\mathbb{Z}_p$.  An element $x\in H^*(G,k)$ is \textit{essential} if it restricts to zero on every proper subgroup $H\subset G$, that is, if $\textnormal{res}^G_Hx=0$ for all $H\subset G$.  The set of all essential elements will be denoted $Ess(G)$, and this is an ideal of $H^*(G,k)$ by the following lemma.\\
\\
\textbf{Lemma 2.1.}  \textit{Let $\mathfrak{M}$ denote the set of maximal subgroups of $G$.  Then $Ess(G)=\textnormal{Ker}\lbrace \textnormal{res}:H^*(G,k)\rightarrow \prod_{M\in\mathfrak{M}}H^*(M,k)\rbrace$}.
\begin{proof}If $u\in Ess(G)$ then $\textnormal{res}^G_Mu=0$ for all $M\in\mathfrak{M}$ by definition of an essential element, so $u$ is contained in the kernel of every restriction $\textnormal{res}^G_M$ and hence is contained in $\textnormal{Ker}\lbrace \textnormal{res}:H^*(G,k)\rightarrow \prod_{M\in\mathfrak{M}}H^*(M,k)\rbrace$.\\
Conversely, if $u\in\textnormal{Ker}\lbrace \textnormal{res}:H^*(G,k)\rightarrow \prod_{M\in\mathfrak{M}}H^*(M,k)\rbrace$ then in particular $\textnormal{res}^G_Mu=0$ for all maximal subgroups $M\subset G$.  Now any proper subgroup $P\subset G$ is contained in some $M$, so $\textnormal{res}^G_Pu=\textnormal{res}^M_P\textnormal{res}^G_Mu=\textnormal{res}^M_P0=0$ for all $P$ and hence $u\in Ess(G)$.
\end{proof}
\indent A theorem of Quillen[Qu] states that if $u\in H^*(G,\mathbb{Z}_p)$ restricts to zero on every elementary abelian $p$-subgroup of $G$, then $u$ is nilpotent.  Quillen's result implies that if $G$ is not elementary abelian then $Ess(G)$ is nilpotent.  However, if $G$ is elementary abelian then it is a fact that the product of the Bocksteins of all nonzero elements of $H^1(G,k)$ is a non-nilpotent essential class.
\\
\\
\textbf{Proposition 2.1.}  \textit{If $G$ is not a $p$-group then $Ess(G)=0$.}
\begin{proof}
For a Sylow $p$-subgroup $P\subset G$ we have $|G:P|$ invertible in $k$ and hence $\textnormal{res}^G_P$ is an injection by Proposition III.10.4[Br] because it maps $H^*(G,k)$ isomorphically onto the set of $G$-invariants in $H^*(P,k)$.  If $u\in Ess(G)$ then $\textnormal{res}^G_Pu=0$ and hence $u=0$ by injectivity.
\end{proof}
\noindent
\textbf{Proposition 2.2.}  \textit{For $p=2$, $Ess(G)=\bigcap \lbrace (x)\,|\, x\in H^1(G,\mathbb{Z}_2)\rbrace$.}
\begin{proof}
A proposition of Marx [Ma, Proposition 2.1] shows that $\textnormal{Ker}(\textnormal{res}^G_H)$ is the principal ideal $(x)$, where $|G:H|=2$ and $x\in H^1(G,\mathbb{Z}_2)$ is a homomorphism $x:G\rightarrow \mathbb{Z}_2$ such that $\textnormal{Ker}x=H$.  Since the maximal subgroups of a $p$-group are the subgroups of index $p$, we see that every nontrivial element $x$ corresponds to some maximal subgroup $M\subset G$ (which has index $2$) with $\textnormal{Ker}x=M$.
\end{proof}

\section{Calculations: Class A}
\textbf{Theorem 3.1.}  \textit{Let $G=\mathbb{Z}_{p^n}$ where $n\ge 2$.  Then $Ess(G)=(x)$.}
\begin{proof}
Its cohomology ring is given by $H^*(\mathbb{Z}_{p^n},k)\cong \Lambda_k[x]\otimes_k k[y]$, with $|x|=1$ and $|y|=2$.\\
\indent Since there is a unique maximal subgroup $H=\mathbb{Z}_{p^{n-1}}$ of $G$, $Ess(G)=\textnormal{Ker}(\textnormal{res}^G_H)$.  As $G$ is not elementary abelian, $Ess(G)$ is nilpotent; thus $y\notin Ess(G)$.  Note that we could also deduce this by viewing $y\in H^2(G,k)$ as a group extension and showing that it restricts to a non-split extension.\\
\indent For $(p=2,\,n=2)$ the restriction map is $\textnormal{res}^G_H:\Lambda_k[x]\otimes_k k[y]\rightarrow k[w]$, with $|w|=1$.  Non-essentiality of $y\in H^*(G,k)$ and dimension considerations force $\textnormal{res}^G_H(y)=w^2$.  In all other cases $(p,n)$ the restriction map is $\textnormal{res}^G_H:\Lambda_k[x]\otimes_k k[y]\rightarrow\Lambda_k[w_1]\otimes_k k[w_2]$, with $|w_1|=1$ and $|w_2|=2$.  Non-essentiality of $y\in H^*(G,k)$ and dimension considerations force $\textnormal{res}^G_H(y)=1\otimes w_2$, noting that $w_1^2=0\in\Lambda_k[w_1]$.\\
\indent In particular, $\textnormal{res}^G_H(y^i)$ is either $w^{2i}$ or $1\otimes w_2^i$ for all $i\in\mathbb{N}$, which are both nonzero elements.  Thus $y^i\notin Ess(G)$ for all $i\in\mathbb{N}$.\\
\indent We can view $x\in H^1(G,\mathbb{Z}_p)$ as a nontrivial homomorphism $x:G\rightarrow \mathbb{Z}_p$ with kernel $H$.  As the restriction map is induced from the inclusion $H\hookrightarrow G$, we have $\textnormal{res}^G_H(x)=0$.  Therefore, $(x)\subseteq Ess(G)$.\\
\indent I claim that $Ess(G)=(x)$.  Indeed, it suffices to show that the induced map $res:H^*(G,k)/(x)\cong k[y]\rightarrow H^*(H,k)$ is injective.  But this is immediate, because it is determined by $res(y)$, which is nontrivial (as explained above).
\end{proof}
\noindent
\textbf{Theorem 3.2.}  \textit{Let $G=\mathbb{Z}_2$.  Then $Ess(G)=(x)$.}
\begin{proof}
Its cohomology ring is given by $H^*(\mathbb{Z}_{2},k)\cong k[x]$, with $|x|=1$.\\
\indent Since the only proper subgroup of $G$ is $\lbrace 1\rbrace$, and all nonzero-degree elements restrict to zero on the trivial group, we have $Ess(G)=(x)$.  This also follows from Proposition 2.2 because there is only a single generating class in $H^1(G,\mathbb{Z}_2)$.
\end{proof}
\noindent
\textbf{Theorem 3.3.}  \textit{Let $G=\mathbb{Z}_p$ where $p>2$.  Then $Ess(G)=(x,y)$.}  
\begin{proof}
Its cohomology ring is given by $H^*(\mathbb{Z}_{p},k)\cong \Lambda_k[x]\otimes_k k[y]$, with $|x|=1$ and $|y|=2$.\\
\indent Since the only proper subgroup of $G$ is $\lbrace 1\rbrace$, and all nonzero-degree elements restrict to zero on the trivial group, we have $Ess(G)=(x,y)$.
\end{proof}

\section{Calculations:  Class B}
\textbf{Theorem 4.1.}  \textit{Let $G=\mathbb{Z}_2\times\mathbb{Z}_2$.  Then $Ess(G)=(x_1^2x_2+x_1x_2^2)$.}
\begin{proof}
Its cohomology ring is given by $H^*(G,k)\cong k[x_1,x_2]$, with $|x_1|=|x_2|=1$.\\
\indent  By Proposition 2.2, $Ess(G)=(x_1)\cap(x_2)\cap(x_1+x_2)=(x_1x_2(x_1+x_2))$.  Alternatively, Lemma 2.2[Ak] states that $Ess(G)$ is generated by $L_2(x_1,x_2)\equiv x_1x_2^2-x_2x_1^2$.
\end{proof}
\noindent
\textbf{Theorem 4.2.}  \textit{Let $G=\mathbb{Z}_p\times\mathbb{Z}_p$, where $p>2$.  Then $Ess(G)=(x_1x_2,x_1y_2-x_2y_1,x_1y_2^p-x_2y_1^p,y_1^py_2-y_1y_2^p)$.}
\begin{proof}
Its cohomology ring is given by $H^*(\mathbb{Z}_{p}\times\mathbb{Z}_p,k)\cong\Lambda_k[x_1,x_2]\otimes_k k[y_1,y_2]$, with $|x_1|=|x_2|=1$ and $|y_1|=|y_2|=2$.\\
\indent Theorem 1.1[Ak] states that $Ess(G)$ is the Steenrod closure of $\Lambda_k^2(G^*)$, where $G^*$ is the dual space of $G$.  The product $x_1x_2$ is a basis for $\Lambda_k^2(G^*)$, so $Ess(G)=(x_1x_2,\beta(x_1x_2),\mathcal{P}^1\beta(x_1x_2),\beta\mathcal{P}^1\beta(x_1x_2))$ which is the Steenrod closure of $x_1x_2$.  Note that $y_1=\beta(x_1)$ and $y_2=\beta(x_2)$ where $\beta$ is the mod-$p$ bockstein homomorphism.  The Steenrod power $\mathcal{P}^1$ sends $y_i$ to $y_i^p$, sends $x_i$ to $0$, and obeys the Cartan formula $\mathcal{P}^1(ab)=\mathcal{P}^1(a)b+a\mathcal{P}^1(b)$.\\
\indent Now $\beta(x_1x_2)=\beta(x_1)x_2+(-1)^{|x_1|}x_1\beta(x_2)=y_1x_2-x_1y_2$, and $\mathcal{P}^1\beta(x_1x_2)=\mathcal{P}^1(x_2)y_1+x_2\mathcal{P}^1(y_1)-\mathcal{P}^1(x_1)y_2-x_1\mathcal{P}^1(y_2)=0+x_2y_1^p-0-x_1y_2^p=x_2y_1^p-x_1y_2^p$.  The result follows.
\end{proof}
\noindent
\textbf{Theorem 4.3.}  \textit{Let $G=\mathbb{Z}_{p^n}\times\mathbb{Z}_p$, where $p>2$ and $n\ge 2$.  Then $Ess(G)=(x_1x_2,x_1y_2)$.}
\begin{proof}
Its cohomology ring is given by $H^*(\mathbb{Z}_{p^n}\times\mathbb{Z}_p,k)\cong\Lambda_k[x_1,x_2]\otimes_k k[y_1,y_2]$, with $|x_1|=|x_2|=1$ and $|y_1|=|y_2|=2$.\\
\indent The maximal subgroups of $G=T\times S=\langle t,s\rangle$ are $K=\langle t^p,s\rangle\cong\mathbb{Z}_{p^{n-1}}\times \mathbb{Z}_p$ and $p$ distinct cyclic groups $\langle ts^i\rangle\cong\mathbb{Z}_{p^n}$, $0\le i\le p-1$.  Up to K$\ddot{u}$nneth isomorphism, $\textnormal{res}^G_T=id_T\otimes\textnormal{res}^S_{\lbrace 1\rbrace}$.  From here it is obvious that $Ker(\textnormal{res}^G_T)=(x_2,y_2)$.  Now $\textnormal{res}^G_K=\textnormal{res}^T_{T_0}\otimes id_S$, where $K=T_0\times S$ and $T_0=\langle t^p\rangle$.  From Theorem 3.1 we know that the kernel of $\textnormal{res}^T_{T_0}$ is the principal ideal generated by $x_1\in H^1(T,k)$, and hence $Ker(\textnormal{res}^G_K)=(x_1)$.\\
\indent  It suffices to consider $Ker(\textnormal{res}^G_H)$ where $H=\langle ts^i\rangle$, for $1\le i\le p-1$.  Write $H^*(H,k)\cong \Lambda_k[x]\otimes_kk[y]$.  Considering 1-dimensional cohomology classes as homomorphisms, the generator $x:H\rightarrow k$ of $H^1(H,k)$ is the image of the generator $\tilde{x}:H\rightarrow \mathbb{Z}_{p^2}$ of $H^1(H,\mathbb{Z}_{p^2})$ under the respective map in the long exact cohomology sequence associated to $k\hookrightarrow\mathbb{Z}_{p^2}\twoheadrightarrow k$.  Thus $\beta(x)=0$, where $\beta$ is the bockstein homomorphism.  Similarly, $\beta(x_1)=0$, but note that $\beta(x_2)=y_2$ since $S\cong\mathbb{Z}_p$.\\
\indent Now $x_1:G\rightarrow k$ is given by $t\mapsto 1$ and $s\mapsto 0$, and $x_2:G\rightarrow k$ is given by $t\mapsto 0$, $s\mapsto 1$.  In particular, $x_1(ts^i)=1$ and $x_2(ts^i)=i$, so that $\textnormal{res}^G_Hx_1=x$ and $\textnormal{res}^G_Hx_2=ix$.  Then $(x_1x_2)\subseteq (ix_1-x_2)\subseteq Ker(\textnormal{res}^G_H)$, where we note that $x_1x_2=-x_1(ix_1-x_2)$.  As the bockstein commutes with restriction, $\textnormal{res}^G_Hy_2=\textnormal{res}^G_H\beta(x_2)=\beta(\textnormal{res}^G_Hx_2)=\beta(ix)=0$.\\
\indent Considering 2-dimensional cohomology classes as group extensions, we have the following commutative diagram\\
\indent\indent\indent\xymatrix{
y_1: & k\ar@{^{(}->}[r]\ar@{=}[d]& \mathbb{Z}_{p^{n+1}}\ar@{->>}[r]& T\\
y_1=\textnormal{inf}^G_T(y_1):& k\ar@{^{(}->}[r]\ar@{=}[d]& \mathbb{Z}_{p^{n+1}}\times \mathbb{Z}_p\ar@{->>}[r]^\pi\ar@{->>}[u]& G\ar@{->>}[u]\\
\textnormal{res}^G_Hy_1:& k\ar@{^{(}->}[r]& E\ar@{->>}[r]\ar@{^{(}->}[u]& H\ar@{^{(}->}[u]^\tau}\\
where $E$ is the pullback $\lbrace (v,w)\in (\mathbb{Z}_{p^{n+1}}\times \mathbb{Z}_p)\times H\;|\; \pi(v)=\tau(w)\rbrace$.  It is then clear that $E\cong\mathbb{Z}_{p^n}\times\mathbb{Z}_p$, so $\textnormal{res}^G_Hy_1=0$.  Thus $Ker(\textnormal{res}^G_H)=(ix_1-x_2,y_1,y_2)$.\\
\indent Putting this all together, $Ess(G)=(x_1)\cap (x_2,y_2)\cap[\bigcap_{i=1}^{p-1}(ix_1-x_2,y_1,y_2)]=(x_1x_2,x_1y_2)$.
\end{proof}
\noindent
\textbf{Theorem 4.4.}  \textit{Let $G=\mathbb{Z}_{2^n}\times\mathbb{Z}_2$, where $n\ge 2$.  Then $Ess(G)=(x_1x_2)$.}
\begin{proof}
Its cohomology ring is given by $H^*(\mathbb{Z}_{2^n}\times\mathbb{Z}_2,k)\cong\Lambda_k[x_1]\otimes_k k[x_2,y_1]$, with $|x_1|=|x_2|=1$ and $|y_1|=2$.\\
\indent By Proposition 2.2 we have $Ess(G)=(x_1)\cap(x_2)\cap(x_1+x_2)=(x_1x_2)$.
\end{proof}

\section{Calculations:  Class C}
\textbf{Theorem 5.1.}  \textit{Let $G=\mathbb{Z}_{2^n}\rtimes \mathbb{Z}_2=\langle t,s\;|\;s^2=t^{2^n}=1,\,sts=t^{2^{n-1}+1}\rangle$ where $n\ge 2$.  Then $Ess(G)=(ab)$.}
\begin{proof}
Its cohomology ring is given by $H^*(\mathbb{Z}_{2^n}\rtimes \mathbb{Z}_2,k)\cong  k[a,b,v,w]/(a^2,v^2,av,ab^2)$, with $|a|=|b|=1$ and $|v|=3$ and $|w|=4$.\\
\indent By Proposition 2.2 we have $Ess(G)=(a)\cap(b)\cap(a+b)$.  From this and the relations in the cohomology ring it is apparent that terms involving $v$ and $w$ do not lie in $Ess(G)$ unless the non-$v$ and non-$w$ elements in the terms lie in $Ess(G)$.  Similarly, $a$ and $b$ do not lie in $Ess(G)$.\\
\indent But $ab$ lies in this intersection because $ab\in(a)\cap (b)$ and $ab=ab+0=ab+a^2=a(b+a)\in(a+b)$.  Thus $(ab)= Ess(G)$.
\end{proof}
\noindent
\textbf{Theorem 5.2.}  \textit{Let $G=\mathbb{Z}_{p^n}\rtimes \mathbb{Z}_p=\langle t,s\;|\;s^p=t^{p^n}=1,\,sts^{-1}=t^{p^{n-1}+1}\rangle$, where $n\ge2$ and $p>2$.  Then $Ess(G)=(a_1b,\ldots,a_{p-1}b,vb,vy)$.}
\begin{proof}
Its cohomology ring is given by\\
$H^*(\mathbb{Z}_{p^n}\rtimes \mathbb{Z}_p,k)\cong k[a_1,\ldots,a_{p-1},b,y,v,w]/(b^2,v^2,a_ia_j,a_iv,a_iy)$, with $|a_i|=2i-1$ and $|b|=1$ and $|y|=2$ and $|v|=2p-1$ and $|w|=2p$.\\
\indent The maximal subgroups of $G$ are $K=\langle t^p,s\rangle\cong\mathbb{Z}_{p^{n-1}}\times \mathbb{Z}_p$ and $p$ distinct cyclic groups $M_i=\langle ts^i\rangle\cong\mathbb{Z}_{p^n}$ for $0\le i<p$, by Proposition IV.4.4[Br].  Let $T=M_0=\langle t\rangle$ and $S=\langle s\rangle$.\\
\indent Abusing notation, I will write $H^*(T,k)\cong \Lambda_k[x]\otimes_kk[z]\cong H^*(S,k)$.  Let $E^{ij}_r$ denote the Lyndon-Hochschild-Serre spectral sequence of the extension for $G$.  From [Di] we know that $a_i$ corresponds to a generator of $E_\infty^{0,2i-1}=E_3^{0,2i-1}$ [including $a_p\equiv v$] and $w$ corresponds to a generator of $E_\infty^{0,2p}=E_3^{0,2p}$ and $b$ corresponds to a generator of $E_\infty^{1,0}=E_2^{1,0}$ and $y$ corresponds to a generator of $E_\infty^{2,0}=E_2^{2,0}$.  Furthermore, $b$ and $y$ are nontrivial images of the inflation map $\textnormal{inf}^G_S:H^*(S,k)\hookrightarrow H^*(G,k)$ which is an injection because of the splitting $S\rightarrow G$.\\
\indent Now $E_2^{ij}=H^i(S,H^j(T,k))\cong H^i(S,k)\otimes_kH^j(T,k)$, so $E_3^{0,j}=\textnormal{Ker}(d_2^{0,j})\subseteq E_2^{0,j}=H^j(T,k)\cong k$ and hence we must actually have $E_3^{0,j}=E_2^{0,j}$.  From this information we see that $a_i$ and $w$ restrict nontrivially on $T$, and $b$ and $y$ restrict nontrivially on $S$.\\
\indent As stated in [Di] we can arrange that $a_i$ and $w$ restrict trivially on $S$.  To see this for $a_i$, consider the quotient map $H^{2i-1}(G,k)\twoheadrightarrow E_\infty^{0,2i-1}=H^{2i-1}(T,k)$.  Let $x_i\in H^{2i-1}(G,k)$ be the element which maps onto the generator $xz^{i-1}\in E_\infty^{0,2i-1}$.  The kernel of this quotient map is $F^1H^{2i-1}(G,k)$, the first filtration submodule, which contains $F^{2i-1}H^{2i-1}(G,k)=H^{2i-1}(S,k)\subseteq F^1H^{2i-1}(G,k)$.  Thus adding an element $\alpha\in H^{2i-1}(S,k)$ to $x_i$ does not have any effect when passing to the quotient.  Then $0=\textnormal{res}^G_S(x_i+\textnormal{inf}^G_S\alpha)=\textnormal{res}^G_Sx_i+\alpha\;\Rightarrow\; \alpha=-\textnormal{res}^G_Sx_i$, and hence we obtain the element $a_i\equiv x_i-\textnormal{inf}^G_S\textnormal{res}^G_Sx_i\in H^{2i-1}(G,k)$ which restricts trivially on $S$ and corresponds to the generator of $E_\infty^{0,2i-1}$.\\
\indent As the composition $T\rightarrow G\rightarrow S$ is the zero map, $\textnormal{res}^G_T\textnormal{inf}^G_{S}=0$ and hence $b$ and $y$ restrict trivially on $T$.\\
\indent Due to dimension considerations, we then must have the following:\\
$\textnormal{res}^G_T(a_i)=z^{i-1}x\;,\;\textnormal{res}^G_T(v)=z^{p-1}x\;,\;\textnormal{res}^G_T(w)=z^p\;,\;\textnormal{res}^G_T(b)=\textnormal{res}^G_T(y)=0$\\
$\textnormal{res}^G_S(a_i)=\textnormal{res}^G_S(v)=\textnormal{res}^G_S(w)=0\;,\; \textnormal{res}^G_S(b)=x\;,\; \textnormal{res}^G_S(y)=z$\\
\indent Note that the image $\textnormal{res}(g)=g'$ of each generator $g$ under the restriction map could actually be a scalar multiple $m_g\cdot g'$ of what is stated, but we have the freedom of forming a new set of generators $\lbrace m_g^{-1}\cdot g\rbrace$ so that $\textnormal{res}(m_g^{-1}\cdot g)=m_g^{-1}[m_g\cdot g']=g'$ and everything else is unaltered.\\
\indent Thus $a_i,b,y^j,v,w^j,a_iw^j,by^j,vw^j\notin Ess(G)$ where $j\in\mathbb{Z}^+$, and these monomial terms map to distinct nonzero [linearly independent] elements (under the direct product of the two restrictions $\textnormal{res}^G_T$ and $\textnormal{res}^G_S$).  In particular, no polynomial involving these monomials could restrict trivially (under the direct product), because no sum of distinct nonzero [linearly independent] elements in $\Lambda_k[x]\otimes_kk[z]$ is the trivial element.\\
\indent Exhausting through all possible combinations of the generators to obtain all monomial terms in $H^*(G,k)$, the ones listed in the previous paragraph are the only ones which do not restrict trivially under either $\textnormal{res}^G_T$ or $\textnormal{res}^G_S$.  So the only elements which might lie in $Ess(G)$ are the polynomials formed by the following monomial terms (with $j,r\in\mathbb{Z}^+$):\\
$X\equiv\lbrace a_ib\;,\; a_ibw^j\;,\; vy^j\;,\; vy^jw^r\;,\; w^jy^r\;,\; bv\;,\; bw^j\;,\; bvy^j\;,\; bvw^j\;,\; bw^jy^r\;,\; bvy^jw^r\rbrace$.
\\
\\
\textbf{Lemma 5.1.}  $\textnormal{res}^G_Ka_i=0$ for all $1\le i\le p$.
\begin{proof}
The $E_2$-page of the Lyndon-Hochschild-Serre spectral sequence for $K=\mathbb{Z}_p\times\mathbb{Z}_{p^{n-1}}=S\times T_0$ is $E_2^{ij}=H^i(S,k)\otimes_kH^j(T_0,k)=E_\infty^{ij}$, where the latter equality is seen to be true from the K$\ddot{u}$nneth isomorphism $H^*(K,k)\cong H^*(S,k)\otimes_kH^*(T_0,k)=\Lambda_k[x_t,x_s]\otimes_kk[z_t,z_s]$.\\
\indent Let $E$ denote the spectral sequence for $G$ and let $\bar{E}$ denote the spectral sequence for $K$.  Let the generators in cohomology denote the corresponding generators in the spectral sequence, and let $Res$ denote the restriction map ($G$ to $K$) at the spectral sequence level.\\
\indent On the $E_2$-page the restriction map in bidegree $(i,0)$ is the identity map $E_2^{i,0}=H^i(S,k)\otimes_kk\rightarrow H^i(S,k)\otimes_kk=\bar{E}_2^{i,0}$.  In particular, $Res(b)=x_s$ and $Res(y)=z_s$.\\
\indent On the $E_2$-page the restriction map in bidegree $(0,j)$ is the tensored map $id_S\otimes \textnormal{res}^T_{T_0}:k\otimes_kH^j(T,k)\rightarrow k\otimes_kH^j(T_0,k)$.  From Theorem 3.1 we know that the kernel of $\textnormal{res}^T_{T_0}$ is the principal ideal generated by $x\in H^1(T,k)$.  In particular, $Res(a_i)=0$ and $Res(w)=z_t^p$.\\
\indent I claim that $\textnormal{res}^G_Ka_i=0$ for all $1\le i\le p$.  To see this, first consider the commutative diagram\\
\indent\indent\indent\indent\xymatrix{
H^{2i-1}(G,k)\ar[d]^{\textnormal{res}^G_K}\ar@{->>}[r] & E_\infty^{0,2i-1}\ar[d]^{Res}\\
H^{2i-1}(K,k)\ar@{->>}[r] & \bar{E}_\infty^{0,2i-1}}\\
Then $\textnormal{res}^G_Ka_i\in F^1H^{2i-1}(K,k)=H^{2i-1}(K,k)/[H^0(S,k)\otimes_kH^{2i-1}(T_0,k)]$ by commutativity of the diagram coupled with $Res(a_i)=0$.  Here the filtration is defined by $F^m=F^{m-1}/[H^{m-1}(S,k)\otimes_kH^{2i-m}(T_0,k)]$.  It suffices to show that $\textnormal{res}^G_Ka_i\in H^0(S,k)\otimes_kH^{2i-1}(T_0,k)$.\\
\indent As stated in [Di] there is an automorphism $\varphi\in Aut(G)$ of order $p-1$ which acts trivially on $S$ but nontrivially on $T$.  Furthermore, this automorphism induces multiplication by $v^i$ on $E_\infty^{0,2i-1}$, where $v$ is a generator of $k^*$, and it acts trivially on $E_\infty^{*,0}$.  Consider the commutative diagram\\
\indent\indent\indent\indent\xymatrix{
H^{2i-1}(G,k)\ar[r]^{\textnormal{res}^G_K}\ar[d]^\varphi & H^{2i-1}(K,k)\ar[d]^\varphi\\
H^{2i-1}(G,k)\ar[r]^{\textnormal{res}^G_K} & H^{2i-1}(K,k)}\\
As stated in [Di] we can further choose each $a_i$ so that $\varphi(a_i)=v^ia_i$.  Then $\varphi(\textnormal{res}^G_Ka_i)=v^i\textnormal{res}^G_Ka_i$, so that $\textnormal{res}^G_Ka_i$ lies in the $v^i$-eigenspace $V_i$ of $H^*(K,k)$.  It suffices to show that $\varphi$ induces multiplication by $v^i$ on $\bar{E}_\infty^{0,2i-1}$, for then $V_i=H^0(S,k)\otimes_kH^{2i-1}(T_0,k)$ and hence $\textnormal{res}^G_Ka_i\in H^0(S,k)\otimes_kH^{2i-1}(T_0,k)$.\\
\indent On $E_\infty^{0,*}=H^*(T,k)\cong \Lambda_k[x]\otimes_kk[z]$ we have $\varphi(x)=vx$ and $\varphi(z)=vz$, so that $\varphi(a_i)=\varphi(z^{i-1}x)=v^{2i-1}z^{i-1}x=v^{2i-1}a_i$.  Here $v$ is the image of $\tilde{v}\in\mathbb{Z}_{p^n}^*=Aut(T)$ under the mod-$p$ restriction $\mathbb{Z}_{p^n}^*\rightarrow\mathbb{Z}_p^*=k^*$, where $\tilde{v}$ is a generator of the unique subgroup of $\mathbb{Z}_{p^n}^*$ of order $p-1$ representing $\varphi|_T$.  We can restrict $\varphi|_T$ to $\varphi|_{T_0}$, sending $\tilde{v}$ to $\tilde{\tilde{v}}$ in the unique subgroup of $\mathbb{Z}_{p^{n-1}}^*=Aut(T_0)$ of order $p-1$, and $\tilde{\tilde{v}}$ maps to $v$ under the mod-$p$ restriction.  Thus $\varphi$ also induces multiplication by $v^i$ on $\bar{E}_\infty^{0,2i-1}$.
\end{proof}
\indent It is now apparent that $I=(a_1b,\ldots,a_{p-1}b,vb,vy)$ \textit{might} lie in $Ess(G)$, but all other monomial terms in $X$ (call that collection $X'$) do not lie in the essential ideal.  These elements of $X'$ map to distinct nonzero [linearly independent] elements under $Res$ on the spectral sequence level (given in Lemma 5.1).  In particular, no polynomial involving the elements of $X'$ could restrict trivially under $Res$, because no sum of distinct nonzero [linearly independent] elements in $\bar{E}_\infty^{ij}=H^i(S,k)\otimes_kH^j(T_0,k)$ is the trivial element.\\
\indent It suffices to compute $\textnormal{res}^G_{M_i}$ on $I$ for $0<i<p$.  Let $H^*(M_i,k)=\Lambda_k[\alpha]\otimes_kk[\beta]$.  Viewing $b\in H^1(G,k)$ as a homomorphism $G\rightarrow k$, we have $b(t)=0$ and $b(s)=1$.  Then $b(ts^i)=b(t)+i\cdot b(s)=i$, so $\textnormal{res}^G_{M_i}b=i\alpha$.  Similarly, $\textnormal{res}^G_{M_i}a_1=\alpha$.  Via dimension considerations we have $\textnormal{res}^G_{M_i}a_j=c_j\beta^{j-1}\alpha$ and $\textnormal{res}^G_{M_i}y=c_y\beta$, where the constants $c_j,c_y\in k$ might possibly depend on $i$.  Then $\textnormal{res}^G_{M_i}(a_jb)=\textnormal{res}^G_{M_i}(vb)=0$ for all $i$.  Furthermore, since $0=a_1y$ we have $0=\textnormal{res}^G_{M_i}(a_1y)=\alpha\cdot c_y\beta$ and hence $c_y=0$, i.e. $\textnormal{res}^G_{M_i}y=0$.  Thus $\textnormal{res}^G_{M_i}(I)=0$ and $Ess(G)=I$.
\end{proof}

\section{Calculations:  Class D}
\textbf{Theorem 6.1.}  \textit{Let $G=D_{2^n}$, where $n\ge3$.  Then $Ess(G)=0$.}
\begin{proof}
Its cohomology ring is given by $H^*(D_{2^n},k)\cong k[x,y,z]/(xy)$, with $|x|=|y|=1$ and $|z|=2$.\\
\indent There are no nontrivial nilpotent elements in $k[x,y,z]/(xy)$.  But $G$ is not elementary abelian, so $Ess(G)$ is nilpotent; thus $Ess(G)=0$.
\end{proof}

\section{Calculations:  Class E}
\textbf{Theorem 7.1.}  \textit{Let $G=Q_8=\langle i,j\;|\; i^4=1,\,iji=j,\,i^2=j^2\rangle$.  Then $Ess(G)=(x^2,y^2)$.}
\begin{proof}
Its cohomology ring is given by $H^*(Q_8,k)\cong k[x,y,z]/(x^2+xy+y^2,x^2y+xy^2)$ with $|x|=|y|=1$ and $|z|=4$.\\
\indent There are three maximal subgroups, $\langle i\rangle$, $\langle j\rangle$, and $\langle ij\rangle$, each isomorphic to $\mathbb{Z}_4$.  We write the cohomology ring of each of these subgroups as $\Lambda_k[w_1]\otimes_kk[w_2]$, where $|w_1|=1$ and $|w_2|=2$.\\
\indent Since $Q_8$ is not elementary abelian, $Ess(G)$ is nilpotent.  Thus $z^i\notin Ess(G)$ for all $i\in \mathbb{N}$.  Now $z\in H^4(G,k)$ is a generator, and the cohomology of the quaternion group is periodic of period $4$ (see [Ei], pg253-254).  Thus $z$ is isomorphic to the generator in $H^0(G,k)$ which doesn't restrict to zero on any proper subgroup (the restriction map is the identity).  Alternatively, since the Tate cohomology $\widehat{H}^*(G,k)$ is a ring and is periodic, the generator $z$ is invertible.  As any restriction map is a ring homomorphism, it must send invertible elements to invertible elements, and so in particular it must send $z$ to a nonzero element.  Due to dimension considerations, we must have $\textnormal{res}(z)=w_2^2$ on all three maximal subgroups.\\
\indent We can view $x,y\in H^1(G,\mathbb{Z}_2)$ as nontrivial homomorphisms $G\rightarrow \mathbb{Z}_2$.  In particular, $x:Q_8\rightarrow \mathbb{Z}_2$ is given by $i\mapsto 1$, $j\mapsto 0$, $ij\mapsto 1$.  Thus $\langle i\rangle\nsubseteq\textnormal{Ker}x$ and $\langle ij\rangle\nsubseteq\textnormal{Ker}x$ and $\langle j\rangle\subseteq\textnormal{Ker}x$.  As the restriction map is induced by the inclusion, we have $\textnormal{res}^G_{\langle i\rangle}(x)=w_1$ and $\textnormal{res}^G_{\langle j\rangle}(x)=0$ and $\textnormal{res}^G_{\langle ij\rangle}(x)=w_1$.\\
\indent Similarly, $y:Q_8\rightarrow \mathbb{Z}_2$ is given by $i\mapsto 0$, $j\mapsto 1$, $ij\mapsto 1$.  Thus $\langle i\rangle\subseteq\textnormal{Ker}y$ and $\langle ij\rangle\nsubseteq\textnormal{Ker}y$ and $\langle j\rangle\nsubseteq\textnormal{Ker}y$.  As the restriction map is induced by the inclusion, we have $\textnormal{res}^G_{\langle i\rangle}(y)=0$ and $\textnormal{res}^G_{\langle j\rangle}(y)=w_1$ and $\textnormal{res}^G_{\langle ij\rangle}(y)=w_1$.\\
\indent In all three restriction maps, $x^2$ and $y^2$ (and hence $xy$) map to $0$ (either $0^2=0$ or $w_1^2=0$).  Thus $x^2,y^2\in Ess(G)$ and $x,y\notin Ess(G)$.  With the restriction maps sending $z$ to $w_2^2$, we have $\textnormal{res}(xz)=w_1w_2^2\ne0$ and $\textnormal{res}(yz)=w_1w_w^2\ne0$.\\
\indent Noting that $x^3=y^3=0$ [indeed, $x^3=x^3+(x^2y+xy^2)=x(x^2+xy+y^2)=x\cdot0=0$], the above calculations imply that the only monomial terms which lie in $Ess(G)$ are $x^2$ and $y^2$ (hence also $xy$, since $xy=x^2+y^2$).\\
\indent I claim that $Ess(G)=(x^2,y^2)$.  Indeed, it suffices to show that the induced map $res:H^*(Q_8,k)/(x^2,y^2)\cong\Lambda_k[x,y]/(xy)\otimes_kk[z]\rightarrow H^*(\langle i\rangle,k)\times H^*(\langle j\rangle,k)\times H^*(\langle ij\rangle,k)\cong H^*(\mathbb{Z}_4,k)^3$ is injective.  We know that no monomials lie in the kernel of this map, so we may restrict our attention to polynomials (at least two terms).  Since we are working under $k=\mathbb{Z}_2$, all terms in the polynomials must be distinct.  Thus any polynomial in the domain is given by a sum of distinct monomial terms, and $res$ is injective on each of these terms.  Thus the image of any polynomial under $res$ is a sum of distinct monomial terms.  Since no sum of distinct elements in $H^*(\mathbb{Z}_4,k)\cong\Lambda_k[w_1]\otimes_kk[w_2]$ is trivial, $res$ is injective on polynomials.
\end{proof}
\noindent
\textit{Remark:}  We could have simply used Proposition 2.2 to arrive at the same answer.  Indeed, we have $Ess(G)=(x)\cap(y)\cap(x+y)$.  From this it is apparent that $x,y,z^i\notin Ess(G)$ for all $i\in \mathbb{N}$.  But $x^2$ and $y^2$ lie in this intersection because $x^2=x\cdot x\in(x)$ and $x^2=xy+y^2=y\cdot (x+y)\in (x+y)\cap(y)$ [same for $y^2$].
\\
\\
\textbf{Theorem 7.2.}  \textit{Let $G=Q_{2^n}=\langle t,s\;|\; t^{2^{n-2}}=s^2,\, t^{2^{n-1}}=1,\, tst=s\rangle$, where $n\ge 4$.  Then $Ess(G)=(x^3)$.}
\begin{proof}
Its cohomology ring is given by $H^*(Q_{2^n},k)\cong k[x,y,z]/(xy,x^3+y^3)$, with $|x|=|y|=1$ and $|z|=4$.\\
\indent There is a unique cyclic subgroup $C=\mathbb{Z}_{2^{n-1}}=\langle t\rangle$ of index $2$ (pg98[Br]), and there are two other maximal subgroups (generalized quaternion), $H=\langle t^2,s\rangle$ and $K=\langle t^2,ts\rangle$.  Let us write $\Lambda_k[w_1]\otimes_kk[w_2]$ for the cohomology ring of $C$, and $k[a,b,c]/(ab,a^3+b^3)$ for the cohomology rings of $H$ and $K$.\\
\indent Note that $x^4=y^4=0$.  Indeed, $x^4=x^3\cdot x=y^3\cdot x=y^2\cdot yx=y^2\cdot 0=0$ (same for $y^4$).\\
\indent We can view $x,y\in H^1(G,\mathbb{Z}_2)$ as nontrivial homomorphisms $G\rightarrow \mathbb{Z}_2$.  In particular, $x:G\rightarrow \mathbb{Z}_2$ is given by $t\mapsto 1$, $s\mapsto 0$, $ts\mapsto 1$.  Thus $K\nsubseteq\textnormal{Ker}x$ and $C\nsubseteq\textnormal{Ker}x$ and $H\subseteq\textnormal{Ker}x$.  As the restriction map is induced by the inclusion, we have $\textnormal{res}^G_{H}(x)=0$ and $\textnormal{res}^G_{K}(x)\ne0$ and $\textnormal{res}^G_{C}(x)=w_1$.  Similarly, $y:Q_8\rightarrow \mathbb{Z}_2$ is given by $t\mapsto 0$, $s\mapsto 1$, $ts\mapsto 1$.  Thus $K\nsubseteq\textnormal{Ker}y$ and $H\nsubseteq\textnormal{Ker}y$ and $C\subseteq\textnormal{Ker}y$.  As the restriction map is induced by the inclusion, we have $\textnormal{res}^G_{C}(y)=0$ and $\textnormal{res}^G_{K}(y)\ne 0$ and $\textnormal{res}^G_{H}(y)\ne 0$.\\
\indent Since $x$ and $y$ agree on where they send the generators of $K$ (hence $K$), we have $\textnormal{res}^G_K(x)=\textnormal{res}^G_K(y)\ne 0$.  Dimension considerations force this image to be either $a$ or $b$ or $a+b$.  But it cannot be $a$ nor $b$, otherwise $0=\textnormal{res}^G_K(0)=\textnormal{res}^G_K(xy)=\textnormal{res}^G_K(x)\textnormal{res}^G_K(y)=a^2\ne 0$ (same for $b$).  Thus $\textnormal{res}^G_K(x)=\textnormal{res}^G_K(y)=a+b$.  In particular, $\textnormal{res}^G_H(x^3)=0^3=0$ and $\textnormal{res}^G_C(x^3)=w_1^3=0$ and $\textnormal{res}^G_K(x^3)=(a+b)^3=(a^3+b^3)+(ab^2+ba^2)=0$, while $\textnormal{res}^G_K(x^2)=\textnormal{res}^G_K(y^2)=(a+b)^2=a^2+b^2\ne 0$.\\
\indent Therefore, $(x^3)=(y^3)\subseteq Ess(G)$, while $x,y,x^2,y^2\notin Ess(G)$ and $xy=0$.\\
\indent Since $G$ is not elementary abelian, $Ess(G)$ is nilpotent.  Thus $z^i\notin Ess(G)$ for all $i\in \mathbb{N}$.  Now $z\in H^4(G,k)$ is a generator, and the cohomology of the generalized quaternions is periodic of period $4$ (see [Ei], pg253-254).  Thus $z$ is isomorphic to the generator in $H^0(G,k)$ which doesn't restrict to zero on any proper subgroup (the restriction map is the identity).  Alternatively, since the Tate cohomology $\widehat{H}^*(G,k)$ is a ring and is periodic, the generator $z$ is invertible.  As any restriction map is a ring homomorphism, it must send invertible elements to invertible elements, and so in particular it must send $z$ to a nonzero element.  Due to dimension considerations, we must have $\textnormal{res}^G_C(z)=w_2^2$ and $\textnormal{res}^G_K(z)=\textnormal{res}^G_H(z)=c$.  In particular, we see that $xz^i,yz^i,x^2z^i,y^2z^i\notin Ess(G)$.\\
\indent I claim that $Ess(G)=(x^3)$.  Indeed, it suffices to show that the induced map $res:H^*(G,k)/(x^3)\cong k[x,y,z]/(xy,x^3,y^3)\rightarrow H^*(H,k)\times H^*(K,k)\times H^*(C,k)$ is injective.  We know that no monomials lie in the kernel of this map, so we may restrict our attention to polynomials (at least two terms).  Since we are working under $k=\mathbb{Z}_2$, all terms in the polynomials must be distinct.  Thus any polynomial in the domain is given by a sum of distinct monomial terms, and $res$ is injective on each of these terms (only under $\textnormal{res}^G_K$ do we have coinciding images for $x$ and $y$).  Thus the image of any polynomial under $res$ is a sum of distinct monomial terms.  Since no sum of distinct elements in $H^*(C,k)$ or in $H^*(H,k)$ or in $H^*(K,k)$ is trivial, $res$ is injective on polynomials.
\end{proof}
\noindent
\textit{Remark:}  Note that we could have simply used Proposition 2.2 to arrive at the same answer.  Indeed, we have $Ess(G)=(x)\cap(y)\cap(x+y)$.  From this and the relations in the cohomology ring it is apparent that $x,y,x^2,y^2,z^i\notin Ess(G)$ for all $i\in \mathbb{N}$.  But $x^3$ (hence $y^3$) lies in this intersection because $x^3\in(x)$ and $x^3=y^3\in(y)$ and $x^3=x^3+xy=x^2(x+y)\in (x+y)$.

\section{Calculations:  Class F}
\textbf{Theorem 8.1.}  \textit{Let $G=\mathbb{Z}_{2^n}\rtimes \mathbb{Z}_2=\langle t,s\;|\;s^2=t^{2^n}=1,\,sts=t^{2^{n-1}-1}\rangle$, where $n\ge 3$.  Then $Ess(G)=(ab)$.}
\begin{proof}
Its cohomology ring is given by\\
$H^*(\mathbb{Z}_{2^n}\rtimes \mathbb{Z}_2,k)\cong k[a,b,y,v,w]/(ay,av,b^2,a^2+ab,v^2+wab+vyb)$, with $|a|=|b|=1$ and $|y|=2$ and $|v|=3$ and $|w|=4$.\\
\indent By Proposition 2.2 we have $Ess(G)=(a)\cap(b)\cap(a+b)$.  Now $ab$ (hence $a^2$) lies in this intersection because $ab\in(a)\cap(b)$ and $ab=ba+0=ba+b^2=b(a+b)\in(a+b)$.  Thus $(ab)\subseteq Ess(G)$.\\
\indent I claim that $Ess(G)=(ab)$.  Indeed, it suffices to show that $(a)\cap(b)\cap(a+b)=0$ in $H^*(G,k)/(ab)\cong \Lambda_k[a,b]\otimes_kk[y,v,w]/(ay,av,ab,v^2+vyb)$.  But this is immediate, because none of the relations in this quotient relate $(a)$ and $(b)$ together.
\end{proof}

\section{Acknowledgements}
\indent\indent I would like to thank Erg$\ddot{u}$n Yalcin for proposing this computational task to me, and Fatma Altunbulak Aksu for reading and commenting on an earlier version of the paper.  I am most grateful towards my undergraduate advisor Ken Brown, who introduced me to the subject of group cohomology and helped me greatly along the way.

\section{References}
\footnotesize{
[Ad1]  A. Adem \& J. Milgram, \textit{Cohomology of Finite Groups} (2nd Ed.), Springer (2004).\\
\newline
[Ad2]  A. Adem \& D. Karagueuzian, \textit{Essential Cohomology of Finite Groups}, Comment. Math. Helv. 72 (1997).\\
\newline
[Ak]  F. Altunbulak Aksu \& D. Green, \textit{Essential Cohomology for Elementary Abelian $p$-Groups}, Journal of Pure and Applied Algebra 213 (2009).\\
\newline
[Br]  K. Brown, \textit{Cohomology of Groups}, Springer GTM 87 (1982).\\
\newline
[Di]  T. Diethelm, \textit{The Mod-p Cohomology Rings of the Nonabelian Split Metacyclic p-Groups}, Arch. Math. 44 (1985).\\
\newline
[Ei]  S. Eilenberg \& H. Cartan, \textit{Homological Algebra}, Princeton University Press (1956).\\
\newline
[Gr]  D. Green, \textit{The Essential Ideal is a Cohen-Macaulay Module}, Proceedings of the Amer. Math. Soc. 133 (2005).\\
\newline
[Ma]  T. Marx, \textit{The Restriction Map in Cohomology of Finite 2-Groups}, Journal of Pure and Applied Algebra 67 (1990).\\
\newline
[Mi]  P. Minh, \textit{Essential Mod-p Cohomology Classes of p-Groups:  An Upper Bound}, Bull. London Math. Soc. 32 (2000)\\
\newline
[Mi2]  P. Minh, \textit{Essential Cohomology and Extraspecial p-Groups}, Trans. Amer. Math. Soc. 353 (2000)\\
\newline
[Qu]  D. Quillen \& B. Venkov, \textit{Cohomology of Finite Groups and Elementary Abelian Subgroups}, Topology 11 (1972).
}
\end{document}